\documentclass[12pt]{amsart}

\usepackage{amsmath}
\usepackage{amssymb}
\usepackage{array}
\usepackage{amsthm}

\renewcommand{\atop}[2]{%
\genfrac{}{}{0pt}{}{#1}{#2}}

\newtheorem{theorem}{Theorem}
\newtheorem{lemma}{Lemma}
\newtheorem{corollary}{Corollary}

\theoremstyle{definition}

\newtheorem{example}{Example}

\renewcommand{\atop}[2]{%
\genfrac{}{}{0pt}{}{#1}{#2}}

\begin{document}

\title{Rational approximations to values of Bell polynomials at  points involving Euler's constant and zeta values.}

\author{Kh.~Hessami Pilehrood}

\address{\begin{flushleft} Mathematics Department,  Faculty of Basic
Sciences,  Shahrekord University,  \\  P.O. Box 115, Shahrekord, Iran \end{flushleft}}


\author{T.~Hessami Pilehrood}
\address{\begin{flushleft} School of Mathematics, Institute for Research  in Fundamental Sciences
(IPM), P.O.Box 19395-5746, Tehran, Iran \end{flushleft}}

\email{hessamik@ipm.ir,  hessamit@ipm.ir, hessamit@gmail.com}


\date{}


\maketitle

\centerline{\small\it Dedicated to the memory of Alfred van der Poorten}

\begin{abstract}
In this paper, we  present
new explicit simultaneous rational approximations converging sub-exponentially to the values of Bell polynomials
at the points of the form $(\gamma, 1! (2a+1)\zeta(2), 2!\zeta(3), \ldots, (m-1)!(a+1+(-1)^ma)\zeta(m)),$
$m=1,2,\ldots,a,$ $a\in{\mathbb N}.$
\end{abstract}


\section{Introduction}

In 2007, A.~I.~Aptekarev and his collaborators \cite{ap,tu} discovered  a sequence of
rational approximations $\tilde{p}_n/\tilde{q}_n$ converging to Euler's constant
$$
\gamma=\lim_{k\to\infty}\left(1+\frac{1}{2}+\cdots+\frac{1}{k}-\log k\right).
$$
sub-exponentially. More precisely,
the numerators $\tilde{p}_n$ and denominators $\tilde{q}_n$ of the approximations are positive integers
generated by the following recurrence relation:
\begin{equation*}
\begin{split}
(16n-15)\tilde{q}_{n+1}&=(128n^3+40n^2-82n-45)\tilde{q}_n \\
&-n^2(256n^3-240n^2+64n-7)\tilde{q}_{n-1}
+n^2(n-1)^2(16n+1)\tilde{q}_{n-2}
\end{split}
\end{equation*}
with the initial conditions
\begin{equation*}
\begin{array}{ccc}
\tilde{p}_0=0, \qquad & \qquad \tilde{p}_1=2, \qquad & \qquad \tilde{p}_2=31, \\
\tilde{q}_0=1,  \qquad &  \qquad \tilde{q}_1=3, \qquad &  \qquad \tilde{q}_2=50
\end{array}
\end{equation*}
and having the following asymptotics:
\begin{align}
\tilde{q}_n&=(2n)!\frac{e^{\sqrt{2n}}}{\sqrt[4]{n}}\left(\frac{1}{\sqrt{\pi}(4e)^{3/8}}+O(n^{-1/2})\right),\nonumber\\[3pt]
\tilde{p}_n-\gamma \tilde{q}_n&=(2n)!\frac{e^{-\sqrt{2n}}}{\sqrt[4]{n}}\left(\frac{2\sqrt{\pi}}{(4e)^{3/8}}+O(n^{-1/2})\right)
\label{eeq01}
\end{align}
The authors \cite{he} found explicit representations for these sequences:
\begin{equation}
\tilde{q}_n=\sum_{k=0}^n\binom{n}{k}^2
\qquad
\tilde{p}_n=\sum_{k=0}^n\binom{n}{k}^2(n+k)!(H_{n+k}+2H_{n-k}-2H_k),
\label{pnqn}
\end{equation}
here $H_n=\sum_{k=1}^n\frac{1}{k}$ is the $n$-th harmonic number, $H_0:=0.$
Formulas (\ref{pnqn}) imply that $\tilde{q}_n$ and $\tilde{p}_n$ are integers divisible by $n!$
and $\frac{n!}{D_n},$ respectively, where $D_n$ denotes the least common
multiple of the numbers $1,2,\ldots, n.$
The linear forms (\ref{eeq01}) do not tend to zero
even after cancelation of their coefficients  by the
 big common factor $\frac{n!}{D_n},$ and therefore this construction
does not allow one to prove the irrationality of $\gamma,$ which is still an open problem.
 Nevertheless,
the linear forms (\ref{eeq01}) present good rational approximations to Euler's constant
$$
\frac{\tilde{p}_n}{\tilde{q}_n}-\gamma=2\pi e^{-2\sqrt{2n}}(1+O(n^{-1/2})) \qquad\mbox{as}
\quad n\to\infty.
$$
In 2009, T.~Rivoal \cite{ri} found another example of rational approximations
to the Euler constant $\gamma,$ viewed as $-\Gamma'(1),$ where $\Gamma$ is the usual Gamma function.
His construction is based on the following third-order recurrence:
\begin{equation*}
\begin{split}
(n&+3)^2(8n+11)(8n+19)y_{n+3}=(n+3)(8n+11)(24n^2+145n+215)y_{n+2}\\
&-(8n+27)(24n^3+105n^2+124n+25)y_{n+1}+(n+2)^2(8n+19)(8n+27)y_n,
\end{split}
\end{equation*}
which provides two sequences of rational numbers $P_n$ and
$Q_n,$ $n\ge 0,$ with the initial values
\begin{equation*}
\begin{array}{lcc}
P_0=-1, \qquad & \qquad P_1=4, \qquad & \qquad P_2=77/4, \\
Q_0=1,  \qquad &  \qquad Q_1=7, \qquad &  \qquad Q_2=65/2
\end{array}
\end{equation*}
such that $\frac{P_n}{Q_n}$ converges to $\gamma.$
The sequences $P_n,$ $Q_n$ satisfy the inclusions (see \cite[Corollary 5]{he1})
$$
n!\,Q_n, \quad n!\,D_n P_n\in {\mathbb Z}
$$
and provide better approximations to $\gamma$
$$
\left|\frac{P_n}{Q_n}-\gamma\right|\le c_0e^{-9/2n^{2/3}+3/2n^{1/3}},
\quad |Q_n|=O(e^{3n^{2/3}-n^{1/3}}) \quad\mbox{as}\quad n\to\infty.
$$
Unfortunately, this convergence
is not fast enough to imply the irrationality of $\gamma.$

In the same paper \cite{ri}, T.~Rivoal considered  a more general construction which is based
on simultaneous Pad\'e approximants to Euler's functions
$$
\mathcal{E}_s(z)=\int_0^{\infty} \frac{\log(t)^{s-1}\,e^{-t}}{z-t}\,dt,
\qquad s\in{\mathbb N}, \,\, z\in {\mathbb C}\setminus [0;+\infty).
$$
This approach allows one in principle to find rational approximations to some other
constants related to higher derivatives of the Gamma function $\Gamma^{(n)}(1),$ $n\ge 2.$
Note that $\Gamma^{(n)}(1)$ can be written as (see \cite[p.~175]{comtet})
$$
\Gamma^{(n)}(1)=Y_n(-\gamma,\zeta(2), -2!\zeta(3), \ldots, (-1)^n (n-1)!\zeta(n)),
$$
where $Y_n$ is the Bell polynomial (see below, for definition and basic properties).
Unfortunately, the details of this construction become rapidly very complicated and
T.~Rivoal presented explicitly rational approximations only for two numbers $\gamma$ and
$\zeta(2)-\gamma^2=\Gamma''(1)-2\Gamma'(1)^2.$ More precisely, he
constructed a sixth order linear recurrence with polynomial
coefficients  of degree $25,$
which has three solutions $\{a_{1,n}\}_{n\ge 0},$ $\{a_{2,n}\}_{n\ge 0},$
and $\{b_{n}\}_{n\ge 0}$ such that $a_{1,n}, a_{2,n}, b_n\in$
$\frac{1}{(3n)!(3n+2)!}\,{\mathbb Z}$ and
$$
\left|\gamma-\frac{a_{1,n}}{b_n}\right|\ll\frac{1}{n^{3/8}b_n},
\qquad \left|\zeta(2)-\gamma^2-\frac{a_{2,n}}{b_n}\right|\ll\frac{1}{n^{3/8}b_n},
$$
$$
|b_n|\sim\frac{c_0}{n^{3/8}}\, \exp(4\sqrt{2} n^{3/4}-5\sqrt{2}/8 n^{1/4})
\qquad\text{as}\quad n\to\infty,
$$
where $c_0$ is some positive constant independent of $n.$
Notice that the better inclusions, namely, $n!^2 b_n, D_n n!^2 a_{1,n}, D_n^2 n!^2 a_{2,n}
\in{\mathbb Z}$ were proved in \cite[Corollary 6]{he1}.

Recently, the authors \cite{he2} gave a new interpretation of Aptekarev's  approximations to
Euler's constant in terms of Meijer $G$-functions and hypergeometric-type series.
This led to finding new rational approximations to $\gamma$ generated
by a second-order inhomogeneous linear
recurrence with polynomial coefficients. The denominators and numerators of these
approximations are given by the formulas
$$
q_n=\sum_{k=0}^n\binom{n}{k}^2 k!, \qquad p_n=\sum_{k=0}^n\binom{n}{k}^2 k! (2H_{n-k}-H_k).
$$
The sequence $\{q_n\}_{n\ge 0}$ satisfies the second-order homogeneous linear recurrence
$$
q_{n+2}-2(n+2)q_{n+1}+(n+1)^2 q_n=0
$$
with $q_0=1,$ $q_1=2,$ and the sequence $\{p_n\}_{n\ge 0}$ is a solution of the second-order
inhomogeneous linear recurrence
$$
q_{n+2}-2(n+2)q_{n+1}+(n+1)^2 q_n=-\frac{n}{n+2}
$$
with $p_0=0,$ $p_1=1.$ Moreover, one has
$$
\frac{p_n}{q_n}-\gamma=e^{-4\sqrt{n}}(2\pi+O(n^{-1/2})) \qquad\text{as}
\quad n\to\infty
$$
and
$$
q_n=n!\,\frac{e^{2\sqrt{n}}}{n^{1/4}}\Bigl(\frac{1}{2\sqrt{\pi e}}+O(n^{-1/2})\Bigr) \qquad\text{as}\quad n\to\infty.
$$
In this paper, we generalize the construction of work \cite{he2} and present explicitly
new simultaneous rational approximations converging sub-exponentially to the values of Bell polynomials
at the points of the form $(\gamma, 1! (2a+1)\zeta(2), 2!\zeta(3), \ldots, (m-1)!(a+1+(-1)^ma)\zeta(m)),$
$m=1,2,\ldots,a,$ $a\in{\mathbb N}.$ Note that our approach is different from that of Rivoal and based
on application of Meijer $G$-functions and complex integrals.

First we recall several known facts related to the Bell polynomials.
The exponential (complete) Bell polynomials (first effectively studied by Eric Temple Bell \cite{bell}
and named in his honor) are the polynomials $Y_n(x_1, \ldots, x_n)$ in an infinite number of variables
$x_1,$ $x_2, \ldots,$ defined by the formal series expansion \cite[\S 3.3]{comtet}:
\begin{equation}
\exp\left(\sum_{m=1}^{\infty}x_m\frac{t^m}{m!}\right)=\sum_{n=0}^{\infty} Y_n(x_1, x_2, \ldots, x_n)
\,\frac{t^n}{n!}.
\label{eq01}
\end{equation}
An explicit representation of $Y_n$ is given by \cite[p.~264]{bell}
$$
Y_n(x_1, \ldots, x_n)=\sum_{\pi(n)}\frac{n!}{k_1! \cdots k_n!}\left(\frac{x_1}{1!}\right)^{k_1}
\cdots \left(\frac{x_n}{n!}\right)^{k_n},
$$
where the summation takes place over all partitions $\pi(n)$ of $n,$ i.e., over all $n$-tuples
of non-negative integers $(k_1, \ldots, k_n)$ such that
$$
\sum_{j=1}^n jk_j=n.
$$
The Bell polynomials satisfy the following recurrence relation (see \cite[p.~263]{bell}):
$$
Y_{n+1}(x_1, \ldots, x_{n+1})=\sum_{k=0}^n\binom{n}{k}x_{k+1}Y_{n-k}(x_1, \ldots, x_{n-k}),
\quad n\ge 0, \qquad Y_0=1,
$$
which implies immediately that the complete Bell polynomials have integer coefficients
and therefore,
\begin{equation}
\frac{n!}{1!^{k_1} k_1! 2!^{k_2} k_2! \cdots n!^{k_n} k_n!}\in {\mathbb Z}
\label{eq00}
\end{equation}
for $k_1, \ldots, k_n$ as above.
The first six complete Bell polynomials are as follows:
$$
Y_0=1, \quad Y_1(x_1)=x_1, \quad  Y_2(x_1, x_2)=x_1^2+x_2, \quad
Y_3(x_1, x_2, x_3)=x_1^3+3x_1x_2+x_3,
$$

$$
Y_4(x_1, x_2, x_3, x_4)=x_1^4+6x_1^2x_2+4x_1x_3+3x_2^2+x_4,
$$

$$
Y_5(x_1, x_2, x_3, x_4, x_5)=x_1^5+10x_1^3x_2+10x_1^2x_3+15x_1x_2^2+5x_1x_4
+10x_2x_3+x_5.
$$
Let $x_1,$ $x_2, \ldots,$ $x_n,\ldots,$  $y_1,$ $y_2, \ldots,$ $y_n,\ldots,$
be two infinite sequences of independent variables. Then from (\ref{eq01}) we get easily
the addition theorem \cite[p.~265]{bell}
$$
Y_n(x_1+y_1, \ldots, x_n+y_n)=[Y(x)+Y(y)]^n.
$$
In ordinary notation this is equivalent to
\begin{equation}
Y_n(x_1+y_1, \ldots, x_n+y_n)=\sum_{k=0}^n\binom{n}{k}Y_k(x_1, \ldots, x_k)Y_{n-k}(y_1, \ldots,
y_{n-k}).
\label{eq02}
\end{equation}
The Bell polynomials play an important role in taking the $n$-th derivative of a composite function.
Namely, the $n$-th derivative of the function $e^{f(x)}$ can be expressed in terms of known quantities
by
\begin{equation}
\left(\frac{d}{dx}\right)^n e^{f(x)}=e^{f(x)}\cdot Y_n(f'(x), f''(x), \ldots, f^{(n)}(x)).
\label{eq03}
\end{equation}
This formula is also known as Fa\`{a} di Bruno's formula for the $n$-th derivative of the composite function.

Let as usual, $(\lambda)_m$ be the  Pochhammer symbol (or the shifted factorial) defined  by $(\lambda)_0=1,$
$(\lambda)_m=\lambda(\lambda+1)\cdots(\lambda+m-1),$ $m\ge 1,$ and $H_n^{(m)}$ be a generalized harmonic number
given by $H_n^{(m)}=\sum_{k=1}^n\frac{1}{k^m}$ and $H_n^{(1)}=H_n.$

Now we can formulate our main result.
\begin{theorem}
Let $a\ge 2$ be an integer. For  $\mu=1, 2, \ldots,a-1$ and any non-negative integer $n,$ define
the following sequences of rational numbers:
\begin{equation}
q_n:=\sum_{k=0}^n\binom{n}{k}^a k!\in {\mathbb Z}, \quad
p_{n,\mu}:=\sum_{k=0}^n\binom{n}{k}^a k! \,Y_{\mu}(r_1(k), r_2(k), \ldots, r_{\mu}(k))\in {\mathbb Q},
\label{eq04}
\end{equation}
where
\begin{equation}
r_m(k):=(m-1)! \left(aH_{n-k}^{(m)}+(-1)^m (a-1)H_k^{(m)}\right),
\qquad k=0,1,\ldots, n.
\label{eq05}
\end{equation}
Let
$$
\alpha_{\mu}:=Y_{\mu}(\gamma, 1!(2a-1)\zeta(2), 2!\zeta(3), \ldots,
(\mu-1)!(a+(-1)^{\mu}(a-1))\zeta(\mu)).
$$
Suppose that the coefficients $b_m(a)$ are defined by the expansion
\begin{equation}
-a\log\Bigl(1+\sum_{m=1}^a\frac{(2-\frac{m+1}{a})_m}{(m+1)!}z^m\Bigr)
-\sum_{m=1}^a\frac{(2-\frac{m}{a})_{m-1}}{m!} z^m
=\sum_{m=1}^a b_m(a) z^m+ O(z^{a+1}), \quad |z|<1.
\label{bma}
\end{equation}
In particular, we have $b_1(a)=-a,$ $b_2(a)=\frac{1-a}{2},$
$b_3(a)=\frac{(1-a)(2a-3)}{6a}.$

Then for every $\mu=1, 2, \ldots, a-1$ there exists a positive constant
$c_{\mu}=c_{\mu}(a)$ such that for any non-negative integer $n$ we have
$$
|p_{n,\mu}-q_n \alpha_{\mu}
|\le
\frac{c_{\mu} n!}{n^{\frac{a}{2}+\frac{1}{2a}}}
\,\exp\left(\sum_{m=1}^{a-1}(-1)^mb_m(a)\cos\Bigl(\frac{2\pi m}{a}\Bigr)n^{1-\frac{m}{a}}\right).
$$
Moreover, $D_n^{\mu}\cdot p_{n,\mu}\in {\mathbb Z}$ and the following asymptotic formula holds:
$$
q_n=\frac{n!}{\sqrt{a} (2\pi)^{\frac{a-1}{2}}\, n^{\frac{a}{2}+\frac{1}{2a}}}
\,\exp\left(\sum_{m=1}^a(-1)^m b_m(a) n^{1-\frac{m}{a}}\right) (1+O(n^{-1/a}))
\quad\text{as}\,\, n\to\infty.
$$
\end{theorem}

The sequences $\{p_{n,\mu}/q_n\}_{n\ge 0},$ $\mu=1,2,\ldots, a-1,$ provide good simultaneous rational
approximations converging sub-exponentially to the numbers $\alpha_{\mu}.$ Note that
$$
\alpha_1=\gamma, \quad \alpha_2=\gamma^2+(2a-1)\zeta(2),\quad
 \alpha_3=\gamma^3+(6a-3)\gamma\zeta(2)+2\zeta(3), \quad
\ldots.
$$
\begin{corollary}
Let $a\ge 2$ be an integer. Then for $\mu=1,2,\ldots, a-1$ we have
\begin{equation*}
\begin{split}
\left|\alpha_{\mu}-\frac{p_{n,\mu}}{q_n}\right|&\le c_{\mu}
\exp\left(\sum_{m=1}^{a-1}(-1)^m b_m(a) \Bigl(\cos(2\pi m/a)-1\Bigr)n^{1-\frac{m}{a}}
\right) \\[3pt]
&<\exp\left(a\Bigl(\cos(2\pi/a)-1\Bigr)n^{1-\frac{1}{a}}(1+o(1))\right),
\end{split}
\end{equation*}
where $c_{\mu}=c_{\mu}(a)$ is a positive constant independent of $n.$
\end{corollary}
In particular, for Euler's constant we have
\begin{corollary}
Let $a\ge 2$ be an integer. Let the sequence $\{q_n\}_{n\ge 0}$ be defined in {\rm (\ref{eq04})}
and $p_{n}=\sum_{k=0}^n\binom{n}{k}^a k!(aH_{n-k}-(a-1)H_k).$ Then
$$
\left|\gamma-\frac{p_{n}}{q_n}\right|<e^{a(\cos(2\pi/a)-1)n^{1-\frac{1}{a}}(1+o(1))}.
$$
\end{corollary}
Consider several examples. For $a=2,$ we get rational approximations to Euler's constant
studied in \cite{he2}.
\begin{example}
For $a=3,$ define three sequences
$$
q_n=\sum_{k=0}^n\binom{n}{k}^3 k!, \qquad
p_{n,1}=\sum_{k=0}^n\binom{n}{k}^3 k!(3H_{n-k}-2H_k),
$$
$$
p_{n,2}=\sum_{k=0}^n\binom{n}{k}^3 k!\left((3H_{n-k}-2H_k)^2+3H_{n-k}^{(2)}+
2H_k^{(2)}\right).
$$
Then we have
$$
\left|\gamma-\frac{p_{n,1}}{q_n}\right|<c_1 e^{-9/2n^{2/3}+3/2n^{1/3}},
\qquad \left|\gamma^2+5\zeta(2)-\frac{p_{n,2}}{q_n}\right|<c_2 e^{-9/2n^{2/3}+3/2n^{1/3}},
$$
and
$$
q_n=\frac{n!}{n^{5/3}}\,e^{3n^{2/3}-n^{1/3}}\left(\frac{e^{1/3}}{2\pi\sqrt{3}}+O(n^{-1/3})\right)
\quad\text{as} \,\, n\to\infty.
$$
Applying Zeilberger's algorithm of creative telescoping (see \cite{PWZ}) it is possible to show
(by the same way as in \cite[Lemma 1]{he2}) that the sequences $p_{n,1},$ $q_n$
are solutions of the third-order homogeneous linear recurrence
\begin{equation*}
\begin{split}
(n+1)(8n-9)f_{n+1}&=(24n^3+13n^2-32n-18)f_n \\
&-n(24n^3-75n^2+52n-5)f_{n-1}
+n(n-1)^3(8n-1)f_{n-2}
\end{split}
\end{equation*}
with the initial conditions
\begin{equation*}
\begin{array}{rrc}
{p}_{0,1}=0, \qquad & \qquad {p}_{1,1}=1, \qquad & \qquad \,\,{p}_{2,1}=13/2, \\
{q}_0=1,  \qquad &  \qquad {q}_1=2, \qquad &  \qquad {q}_2=11
\end{array}
\end{equation*}
and the sequence $p_{n,2}$ is a solution of the third-order inhomogeneous linear recurrence
\begin{equation*}
\begin{split}
(n+1)(8n-9)&f_{n+1}=(24n^3+13n^2-32n-18)f_n
-n(24n^3-75n^2+52n-5)f_{n-1}\\
&+n(n-1)^3(8n-1)f_{n-2}+2(8n^4-17n^3+74n^2-12n-9)/(n(n+1))
\end{split}
\end{equation*}
with the initial values $p_{0,2}=0,$ $p_{1,2}=18,$ $p_{2,2}=95.$
\end{example}
\begin{example}
For $a=4,$ put
$$
q_n=\sum_{k=0}^n\binom{n}{k}^4 k!, \qquad
p_{n,1}=\sum_{k=0}^n\binom{n}{k}^4 k!(4H_{n-k}-3H_k),
$$
$$
p_{n,2}=\sum_{k=0}^n\binom{n}{k}^4 k!(r_1^2(k)+r_2(k)),
\,\quad
p_{n,3}=\sum_{k=0}^n\binom{n}{k}^4 k!(r_1^3(k)+3r_1(k)r_2(k)+r_3(k)),
$$
where $r_m(k)$ is defined in (\ref{eq05}). Then we have
$$
\left|\gamma-\frac{p_{n,1}}{q_n}\right|<c_1 e^{-4n^{3/4}+3n^{1/2}-5/8n^{1/4}},
\quad
\left|\gamma^2+7\zeta(2)-\frac{p_{n,2}}{q_n}\right|<c_2 e^{-4n^{3/4}+3n^{1/2}-5/8n^{1/4}},
$$
$$
\left|\gamma^3+21\gamma\zeta(2)+2\zeta(3)-\frac{p_{n,3}}{q_n}\right|<c_3 e^{-4n^{3/4}+3n^{1/2}-5/8n^{1/4}}
$$
and
$$
q_n=\frac{c_0 n!}{n^{17/8}}\, e^{4n^{3/4}-3/2n^{1/2}+5/8n^{1/4}} (1+O(n^{-1/4})) \quad
\text{as}\,\,\, n\to\infty.
$$
Applying Zeilberger's algorithm of creative telescoping it is easy to show that the sequences
$q_n,$ $p_{n,1},$ $p_{n,2}$ satisfy the fourth-order homogeneous linear recurrence
\begin{equation*}
\begin{split}
&(n+2)^2(729n^4-162n^3-171n^2-4n+6)f_{n+2} \\
&=(2916n^7+14661n^6+20862n^5+947n^4-13008n^3-2370n^2+1320n+312)f_{n+1} \\
&-\!(4374n^8\!-\!18468n^7\!-\!82674n^6\!-\!85776n^5\!-\!13062n^4\!+\!24204n^3\!+\!13528n^2\!+\!2680n\!+\!168)f_n \\
&+n^2(2916n^7+28512n^6+61848n^5+37667n^4-12898n^3-17463n^2-2692n+398)f_{n-1} \\
&-n^2(n-1)^4(729n^4+2754n^3+3717n^2+2084n+398)f_{n-2}
\end{split}
\end{equation*}
with the initial conditions
\begin{align*}
q_0 &=1,  &   q_1 &=2,  &   q_2 &=19, &   q_3 &=250, \\
p_{0,1} &=0,  &  p_{1,1} &=1,  &  p_{2,1} &=13,  &   p_{3,1} &=409/3, \\
p_{0,2} &=0,  &  p_{1,2} &=32, &  p_{2,2} &=217, &   p_{3,2} &=26444/9,
\end{align*}
and the sequence $p_{n,3}$ satisfies the fourth-order inhomogeneous linear recurrence
\begin{equation*}
\begin{split}
&(n+2)^2(729n^4-162n^3-171n^2-4n+6)f_{n+2} \\
&=(2916n^7+14661n^6+20862n^5+947n^4-13008n^3-2370n^2+1320n+312)f_{n+1} \\
&-\!(4374n^8\!-\!18468n^7\!-\!82674n^6\!-\!85776n^5\!-\!13062n^4\!+\!24204n^3\!+\!13528n^2\!+\!2680n\!+\!168)f_n \\
&+n^2(2916n^7+28512n^6+61848n^5+37667n^4-12898n^3-17463n^2-2692n+398)f_{n-1} \\
&-n^2(n-1)^4(729n^4+2754n^3+3717n^2+2084n+398)f_{n-2}\\
&-6(729n^{10}+2754n^9-17424n^8-179680n^7-490669n^6-549106n^5-194460n^4\\
&+100424n^3+105332n^2+30840n+3184)/(n(n+1)^2(n+2))
\end{split}
\end{equation*}
with the initial values $p_{0,3}=0,$ $p_{1,3}=60,$ $p_{2,3}=402,$ $p_{3,3}=50761/9.$
\end{example}

\vspace{0.3cm}

\section{Analytical construction}

\vspace{0.3cm}

Let $a\ge 2$ be an integer.
Let us consider the function
$$
F(n,t)=\frac{n!^a}{\Gamma^{a-1}(t+1)\Gamma^a(n-t+1)}, \qquad n=0,1,2,\ldots,
$$
and for each integer $\mu,$ $0\le\mu\le a-1,$ define
$$
F_{n,\mu}:=F_{n,\mu,a}:=\sum_{k=0}^n\left(\left.\frac{d}{dt}\right)^{\mu} F(n,t)\right|_{t=k}.
$$
\begin{lemma} \label{l1}
Let $a\ge 2$ be an integer. Then for each $\mu=1,2,\ldots, a-1$ there exist $\mu$ constants
$\lambda_{\mu, \nu},$ $\nu=1,\ldots, \mu,$ independent of $n$ such that
for $n=0,1,2,\ldots$ we have
$$
p_{n,\mu}-q_n \alpha_{\mu}
=\sum_{\nu=1}^{\mu} \lambda_{\mu, \nu} F_{n, \nu},
$$
$F_{n,0}=q_n$
and $D_n^{\mu}\cdot p_{n,\mu}\in {\mathbb Z},$ where the sequences $q_n$ and $p_{n,\mu}$ are
defined in {\rm (\ref{eq04})}.
\end{lemma}
\begin{proof}
First let us define the function
$$
f(t):=a\log n!-a\log\Gamma(n+1-t)-(a-1)\log\Gamma(t+1), \qquad 0\le t\le n.
$$
Then we have
\begin{equation}
f'(t)=a\psi(n+1-t)-(a-1)\psi(t+1)
\label{eq06}
\end{equation}
and for $m\ge 2,$
\begin{equation}
f^{(m)}(t)=(-1)^{m-1} a\psi^{(m-1)}(n+1-t)-(a-1)\psi^{(m-1)}(t+1).
\label{eq07}
\end{equation}
Using the well-known formula for the derivatives of $\psi(t),$
\begin{equation}
\psi^{(m)}(t)=(-1)^{m+1} m! \zeta(m+1,t), \qquad m\in {\mathbb N},
\label{eq08}
\end{equation}
where
$$
\zeta(s,t)=\sum_{j=0}^{\infty}\frac{1}{(j+t)^{s}},
$$
is the Hurwitz zeta function, we get
$$
f^{(m)}(t)=(-1)^m (m-1)! \left((-1)^{m-1}a\zeta(m,n+1-t)-(a-1)\zeta(m,t+1)\right), \quad m\ge 2.
$$
We will be interested in the values of $f^{(m)}(t)$ at the integer points $t=k,$  $0\le k\le n.$
Then for $m=1$ taking into account the well-known properties of the digamma function
$$
\psi(1)=-\gamma, \qquad \psi(n+1)=H_n-\gamma, \qquad n\in {\mathbb N},
$$
by (\ref{eq06}), we have
\begin{equation}
f'(k)=-\gamma+aH_{n-k}-(a-1)H_k,
\label{eq09}
\end{equation}
 and from (\ref{eq07}), (\ref{eq08}) for any integer $m\ge 2$ we get
\begin{equation}
f^{(m)}(k)=(m-1)!\Bigl(((-1)^{m-1}(a-1)-a)\zeta(m)+aH_{n-k}^{(m)}+
(-1)^m(a-1)H_k^{(m)}\Bigr).
\label{eq10}
\end{equation}
Now if denote $r_m(k)$ by formula (\ref{eq05}) we get $D_n^m\cdot r_m(k)\in {\mathbb Z},$
and moreover, from (\ref{eq09}), (\ref{eq10}) we obtain
\begin{equation*}
r_m(k)=\begin{cases}
f'(k)+\gamma,     & \text{if} \quad m=1; \\
f^{(m)}(k)+(m-1)! (a+(-1)^m(a-1))\zeta(m), & \text{if} \quad m\ge 2.
\end{cases}
\end{equation*}
Now let us notice that $F(n,t)=e^{f(t)}$ and for calculating the $\mu$-th derivative
of $F(n,t)$ we can apply the Fa\`{a} di Bruno formula (\ref{eq03}) to get:
$$
\left(\frac{d}{dt}\right)^{\mu} F(n,t)=e^{f(t)}\cdot Y_{\mu}(f'(t), f''(t), \ldots, f^{(\mu)}(t)),
\qquad 1\le\mu\le a-1.
$$
Now considering the value $Y_{\mu}(r_1(k), r_2(k), \ldots, r_{\mu}(k))$ and applying the addition formula
(\ref{eq02}), we obtain
\begin{equation*}
\begin{split}
&Y_{\mu}(r_1(k), r_2(k), \ldots, r_{\mu}(k))
=\sum_{\nu=0}^{\mu}\binom{\mu}{\nu}
Y_{\nu}(f'(k), f''(k), \ldots, f^{(\nu)}(k)) \\
&\times Y_{\mu-\nu}(\gamma, 1! (2a-1)\zeta(2), 2! \zeta(3),
\ldots, (\mu-\nu-1)! (a+(-1)^{\mu-\nu}(a-1))\zeta(\mu-\nu)) .
\end{split}
\end{equation*}
This implies that
\begin{equation*}
\begin{split}
&Y_{\mu}(r_1(k), r_2(k), \ldots, r_{\mu}(k))-
Y_{\mu}(\gamma, 1!(2a-1)\zeta(2), \ldots, (\mu-1)! (a+(-1)^{\mu}(a-1))\zeta(\mu)) \\
&\qquad\qquad\qquad\qquad\qquad=\sum_{\nu=1}^{\mu}\binom{\mu}{\nu}Y_{\nu}(f'(k), f''(k), \ldots, f^{(\nu)}(k)) \\
&\times Y_{\mu-\nu}(\gamma, 1! (2a-1)\zeta(2),
\ldots, (\mu-\nu-1)! (a+(-1)^{\mu-\nu}(a-1))\zeta(\mu-\nu)), \quad 0\le k\le n.
\end{split}
\end{equation*}
Now multiplying both sides of the last equality by $e^{f(k)}=k!\binom{n}{k}^a$
and summing over $k=0,1,2, \ldots, n$ we get
\begin{equation*}
\begin{split}
&\qquad\qquad p_{n,\mu}-q_n\cdot Y_{\mu}(\gamma, 1!(2a-1)\zeta(2), \ldots,
(\mu-1)! (a+(-1)^{\mu}(a-1))\zeta(\mu)) \\[3pt]
&=\sum_{\nu=1}^{\mu}\binom{\mu}{\nu}Y_{\mu-\nu}(\gamma, 1! (2a-1)\zeta(2),
\ldots, (\mu-\nu-1)! (a+(-1)^{\mu-\nu}(a-1))\zeta(\mu-\nu)) \\
&\qquad\qquad\qquad\qquad\qquad\qquad\times\sum_{k=0}^n\left.\left(\frac{d}{dt}\right)^{\nu}F(n,t)\right|_{t=k}
\end{split}
\end{equation*}
or
\begin{equation*}
p_{n,\mu}-q_n \alpha_{\mu} 
=\sum_{\nu=1}^{\mu}
\lambda_{\mu,\nu} F_{n,\nu},
\end{equation*}
where $q_n$ and $p_{n,\mu}$ are defined in (\ref{eq04}) and the constants $\lambda_{\mu,\nu}$
given by
$$
\lambda_{\mu,\nu}=\binom{\mu}{\nu}\cdot Y_{\mu-\nu}(\gamma, 1!(2a-1)\zeta(2), \ldots,
(\mu-\nu-1)! (a+(-1)^{\mu-\nu}(a-1))\zeta(\mu-\nu))
$$
are independent of $n.$ To prove the inclusion $D_n^{\mu}\cdot p_{n,\mu}\in {\mathbb Z},$
we consider an arbitrary monomial of $D_n^{\mu}\cdot Y_{\mu}(r_1(k),\ldots,r_m(k)),$
which has the form
\begin{equation*}
\begin{split}
&D_n^{\mu}\, \frac{\mu!}{k_1!\cdots k_{\mu}!}\left(\frac{r_1(k)}{1!}\right)^{k_1}
\cdots\left(\frac{r_{\mu}(k)}{\mu!}\right)^{k_{\mu}}=
\frac{\mu!}{1!^{k_1} k_1!\cdots \mu!^{k_{\mu}} k_{\mu}!} \\[5pt]
&\times(D_n\cdot r_1(k))^{k_1}\cdot (D_n^2\cdot r_2(k))^{k_2}\cdots (D_n^{\mu}\cdot r_{\mu}(k))^{k_{\mu}},
\end{split}
\end{equation*}
where $k_1, \ldots, k_{\mu}$ are non-negative integers such that
$1\cdot k_1+2\cdot k_2+3\cdot k_3+\cdots+\mu\cdot k_{\mu}=\mu.$
Now by (\ref{eq00}) and the fact that $D_n^m\cdot r_m(k)\in {\mathbb Z},$ we get the required inclusion.
\end{proof}

Now for each integer $\mu,$ $0\le\mu\le a-1,$ and $u\in {\mathbb R}$ define a complex integral
\begin{equation}
\begin{split}
I_{n,\mu}(u)&:=\frac{1}{2\pi i}\int_L F(n,t)\left(\frac{\pi}{\sin\pi t}\right)^{\mu+1}e^{i\pi tu}\,dt \\
&=(-1)^{(\mu+1)n} n!^aG_{a,a-1}^{0,\mu+1}\left(
\left.\atop{n+1, \ldots,
n+1}{0, \ldots,  0}\right|e^{i\pi u}
\right),
\label{eq11}
\end{split}
\end{equation}
where $L$ is a loop beginning and ending at $-\infty$ and encircling the points
$n, n-1, n-2, \ldots$ once in the positive direction.
Without loss of generality we can assume that $L$ is located symmetrically with respect to
the real axis.
The integral converges according
to the definition of the Meijer $G$-function (see \cite[\S 5.2]{lu}).
Moreover, if $\mu=a-1$ and $|u|<1/2$ we can also choose the contour of integration
as a vertical line going from $c-i\infty$ to $c+i\infty,$
where $c>n$ is an arbitrary constant.

Let us also put
\begin{equation}
\widetilde{I}_{n,\mu}(u):=\sum_{k=0}^n\,\underset{t=k}{\rm res}\,\left(
F(n,t)\Bigl(\frac{\pi}{\sin\pi t}\Bigr)^{\mu+1} e^{i\pi tu}\right),
\qquad 0\le\mu\le a-1.
\label{tilde}
\end{equation}
\begin{lemma} \label{l2}
Let $a\ge 2,$ $0\le\mu\le a-1$ be integers. Then
$$
\widetilde{I}_{n,\mu}(u)=\begin{cases}
I_{n,\mu}(u),     & \quad\text{if} \quad 0\le\mu\le a-2; \\
I_{n,\mu}(u)+O(n^{-a}), & \quad\text{if} \quad \mu=a-1;
\end{cases}
$$
where the constant in $O$ is absolute.
\end{lemma}
\begin{proof}
First note that
$$
F(n,t)\left(\frac{\pi}{\sin\pi t}\right)^{\mu+1} e^{i\pi tu}=
\frac{(-1)^{na}n!^a}{\pi^{a-\mu-1}} \Gamma(t-n)\left(
\frac{\Gamma(t-n)}{\Gamma(t+1)}\right)^{a-1} (\sin\pi t)^{a-\mu-1}\, e^{i\pi tu}.
$$
Then, since on the segment ${\rm Re}\, t=-N-1/2,$ $|{\rm Im}\, t|\le y_0$
with a sufficiently large integer $N,$ we have
$$
|\Gamma(t-n)|\le |\Gamma({\rm Re}\,(t-n))|=|\Gamma(-N-1/2-n)|=
\frac{\pi}{\Gamma(N+1/2+n)}=O(e^{-N\log N+N}),
$$
$$
\left|\frac{\Gamma(t-n)}{\Gamma(t+1)}\right|=O\left(\frac{1}{N^{n+1}}\right), \qquad
|(\sin\pi t)^{a-\mu-1} e^{i\pi tu}|=O(1),
$$
it follows that the integral $I_{n,\mu}$ can be evaluated as a sum of residues
at the singular points lying inside the loop $L.$ The integrand in (\ref{eq11}) has poles of order
$\mu+1$ at the points $0,1,\ldots,n$ and moreover, if $\mu=a-1,$ it has additional
simple poles at the negative integers. Therefore, for $\mu=a-1,$ we have
\begin{equation*}
\begin{split}
I_{n,a-1}(u)&=\widetilde{I}_{n,a-1}(u)+(-1)^{an}n!^a\sum_{k<0}\,\underset{t=k}{\rm res}\,
\left(\frac{\Gamma(t-n) e^{i\pi tu}}{(t(t-1)\cdots(t-n))^{a-1}}\right) \\
&=\widetilde{I}_{n,a-1}(u)+(-1)^{an}n!^a\sum_{k<0}\,\underset{t=k}{\rm res}\,
\left(\frac{e^{i\pi tu}}{\Gamma(n-t+1)(t(t-1)\cdots(t-n))^{a-1}}\frac{\pi}{\sin\pi t}\right) \\
&=\widetilde{I}_{n,a-1}(u)+n!^a\sum_{k=1}^{\infty}\frac{(-1)^{(u+1)k+a-1}}{(n+k)! (k(k+1)\cdots(k+n))^{a-1}} \\
&=\widetilde{I}_{n,a-1}(u)+\frac{1}{(n+1)^a}\sum_{k=0}^{\infty}\frac{(-1)^{(u+1)k+a-1}k!^{a-1}}{(n+2)_a^k}.
\end{split}
\end{equation*}
Finally, since
$$
\left|\sum_{k=0}^{\infty}\frac{(-1)^{(u+1)k+a-1}k!^{a-1}}{(n+2)_a^k}\right|
\le \sum_{k=0}^{\infty} \frac{k!^{a-1}}{(n+2)_a^k}\le\sum_{k=0}^{\infty}\frac{1}{k!}=e,
$$
we get the desired assertion.
\end{proof}

\section{Bernoulli polynomials}

The generalized Bernoulli polynomials $B_n^{(m)}(x)$ of order $m,$ where $m$ is a positive integer,
are defined by the generating formula (see \cite[\S 2.8]{lu})
\begin{equation}
\frac{z^m  e^{xz}}{(e^z-1)^m}=\sum_{n=0}^{\infty} B_n^{(m)}(x) \frac{z^n}{n!}, \qquad\qquad |z|<2\pi.
\label{eq12}
\end{equation}
Numerous properties of these polynomials can be deduced directly from formula (\ref{eq12}).
We mention here only two of them which will be useful in the sequel. For a detailed study
of Bernoulli polynomials, see \cite{no}. For example, comparing powers of $z$ on both sides
of the equality
$$
\frac{z^m e^{yz}}{(e^z-1)^m}\cdot e^{xz}=\frac{z^m e^{(y+x)z}}{(e^z-1)^m}
$$
leads to the addition formula
\begin{equation}
B_n^{(m)}(x+y)=\sum_{k=0}^n\binom{n}{k}B_k^{(m)}(y) x^{n-k}.
\label{eq13}
\end{equation}
Differentiating both sides of (\ref{eq12}) with respect to $z$ and comparing powers of $z$
we get the recursion formula
\begin{equation}
m B_n^{(m+1)}(x)=(m-n) B_n^{(m)}(x)+n(x-m)B_{n-1}^{(m)}(x).
\label{eq14}
\end{equation}
Setting $n=m$ in (\ref{eq14}) we get $B_m^{(m+1)}(x)=(x-m)B_{m-1}^{(m)}(x),$ which implies
\begin{equation}
B_m^{(m+1)}(x)=(x-1)(x-2)\cdots (x-m).
\label{eq15}
\end{equation}
\begin{lemma} \label{l3}
$1)$ Let $m$ be a positive integer. Then the following series expansion:
$$
\frac{z^m}{\sin^m z}=\sum_{n=0}^{\infty}\frac{(-1)^n 4^n B_{2n}^{(m)}(m/2)}{(2n)!}\, z^{2n},
\qquad\qquad |z|<\pi,
$$
holds.

$2)$ Let $m$ be a positive even integer. Then
\begin{equation}
\sum_{k=0}^m\binom{m}{k}B_{k}^{(m+1)}\left(\frac{m+1}{2}\right) 2^k=0.
\label{eq16}
\end{equation}
\end{lemma}
\begin{proof}
Replacing $z$ by $2z$ in (\ref{eq12}) and using the formula $\sinh z=(e^z-e^{-z})/2$
we get
$$
\frac{z^m e^{(2x-m)z}}{\sinh^m z}=\sum_{n=0}^{\infty}B_n^{(m)}(x) \frac{2^n z^n}{n!}, \qquad\qquad |z|<\pi.
$$
Setting $x=m/2$ in the above equality we have
$$
\frac{z^m}{\sinh^m z}=\sum_{n=0}^{\infty}B_n^{(m)}\left(\frac{m}{2}\right) \frac{2^n z^n}{n!}, \qquad\qquad
|z|<\pi.
$$
Taking into account that $z/\sinh z$ is an even function, we get that for any positive integer $m,$
\begin{equation}
B_{2n+1}^{(m)}\left(\frac{m}{2}\right)=0, \qquad\qquad n=0,1,2,\ldots,
\label{eq17}
\end{equation}
and therefore,
\begin{equation}
\frac{z^m}{\sinh^m z}=\sum_{n=0}^{\infty}B_{2n}^{(m)}\left(\frac{m}{2}\right) \frac{4^n z^{2n}}{(2n)!}, \qquad\qquad
|z|<\pi.
\label{eq18}
\end{equation}
Now replacing $z$ by $iz$ in (\ref{eq18})   and recalling that $\sin z=-i\sinh(iz),$ we get the required expansion.

To prove equality (\ref{eq16}), we consider the addition formula  (\ref{eq13}) with $m$ and $n$
replaced by $m+1$ and $m,$ respectively. Then setting $x=1/2,$ $y=(m+1)/2$ we have
$$
\sum_{k=0}^m\binom{m}{k}B_k^{(m+1)}\left(\frac{m+1}{2}\right) 2^k=2^m B_m^{(m+1)}\left(
\frac{m+2}{2}\right).
$$
Now taking into account that $m$ is even and applying  (\ref{eq15}), we get $B_m^{(m+1)}(m/2+1)=0,$
and the lemma is proved.
\end{proof}

\section{Properties of the integrals $I_{n,\mu}(u)$}

\vspace{0.3cm}

\begin{lemma} \label{l4}
Let $a\ge 2,$ $0\le\mu\le a-1$ be integers. Then we have
\begin{equation*}
F_{n,\mu}=\begin{cases}
\widetilde{I}_{n,0}(1),     & \quad\text{if} \quad \mu=0; \\[3pt]
\sum\limits_{j=0}^{[\mu/2]} c_{2j+1,\mu}\widetilde{I}_{n,2j+1}(0), & \quad\text{if} \quad
\mu \,\,\, \text{is odd}; \\
\sum\limits_{j=1}^{\mu/2} c_{2j,\mu}{\rm Re}\,\widetilde{I}_{n,2j}(1), & \quad\text{if} \quad
\mu\ge 2 \,\,\,  \text{is even};
\end{cases}
\end{equation*}
where $c_{\mu, \mu}\ne 0$ and $c_{j,\mu},$ $0\le j\le \mu,$ are real constants independent of $n.$
\end{lemma}
\begin{proof}
The proof is by induction on $\mu.$ For $\mu=0,$ we deduce easily from (\ref{tilde}):
\begin{equation*}
\widetilde{I}_{n,0}(1)=\sum_{k=0}^n\underset{t=k}{\rm res}\,\left(F(n,t) \, e^{i\pi t}
\cdot\frac{\pi}{\sin\pi t}\right)=\sum_{k=0}^n F(n,k)=
\sum_{k=0}^n\binom{n}{k}^a k!=F_{n,0}.
\end{equation*}
Similarly, for $\mu=1,$ we have
$$
\widetilde{I}_{n,1}(0)=\sum_{k=0}^n\underset{t=k}{\rm res}\,\left(
F(n,t)\cdot\Bigl(\frac{\pi}{\sin\pi t}\Bigr)^2\right)=
\sum_{k=0}^n\left.\frac{d}{dt}F(n,t)\right|_{t=k}=F_{n,1}.
$$
For $\mu=2,$ by (\ref{tilde}), we get
$$
\widetilde{I}_{n,2}(1)=\sum_{k=0}^n\underset{t=k}{\rm res}\,\left(F(n,t)\, e^{i\pi t}
\cdot\Bigl(\frac{\pi}{\sin\pi t}\Bigr)^3\right).
$$
Then using the expansions (see Lemma \ref{l3})
\begin{equation}
\left(\frac{\pi}{\sin\pi t}\right)^{\mu+1}=(-1)^{(\mu+1)k}
\sum_{j=0}^{[\mu/2]}\frac{(-1)^j (2\pi)^{2j} B_{2j}^{(\mu+1)}\bigl(\frac{\mu+1}{2}\bigr)}%
{(2j)! \,(t-k)^{\mu+1-2j}}+O(1),
\label{eq19}
\end{equation}
\begin{equation}
e^{i\pi t}=e^{i\pi k} e^{i\pi(t-k)}=(-1)^k\sum_{j=0}^{\mu}\frac{(i\pi)^j}{j!} (t-k)^j
+O\bigl((t-k)^{\mu+1}\bigr),
\label{eq20}
\end{equation}
\begin{equation}
F(n,t)=\sum_{j=0}^{\mu}\frac{d^j}{dt^j}\left. F(n,t)\right|_{t=k} \frac{(t-k)^j}{j!}+
O\bigl((t-k)^{\mu+1}\bigr)
\label{eq21}
\end{equation}
in a neighborhood of the integer point $t=k$ we obtain
$$
\underset{t=k}{\rm res}\,\left(F(n,t)\, e^{i\pi t}
\cdot\Bigl(\frac{\pi}{\sin\pi t}\Bigr)^3\right)=\frac{1}{2}\frac{d^2}{dt^2}\left.F(n,t)\right|_{t=k}+
i\pi\frac{d}{dt} \left.F(n,t)\right|_{t=k}
$$
and therefore, $F_{n,2}=2 {\rm Re}\, \widetilde{I}_{n,2}(1).$
Now assume that $\mu>2$ and the formula holds for $0,1,2,\ldots, \mu-1.$ We will prove it for
$\mu.$ If $\mu$ is odd we have
$$
\widetilde{I}_{n,\mu}(0)=\sum_{k=0}^n\underset{t=k}{\rm res}\,\left(
F(n,t)\cdot\Bigl(\frac{\pi}{\sin\pi t}\Bigr)^{\mu+1}\right)
$$
and from (\ref{eq19}), (\ref{eq21}) we deduce that
\begin{equation*}
\widetilde{I}_{n,\mu}(0)=\frac{1}{\mu!}F_{n,\mu}+\sum_{j=1}^{[\mu/2]} d_{2j,\mu} F_{n, \mu-2j}.
\end{equation*}
Hence, by the induction hypothesis, we conclude that
$$
F_{n,\mu}=\mu!\widetilde{I}_{n,\mu}(0)-\mu!\sum_{j=1}^{[\mu/2]} d_{2j,\mu} F_{n, \mu-2j}
= \mu!\widetilde{I}_{n,\mu}(0)+\sum_{l=0}^{[\mu/2]-1}c_{2l+1,\mu}\widetilde{I}_{n,2l+1}(0),
$$
as required. If $\mu$ is even, then  by (\ref{tilde}), we have
$$
\widetilde{I}_{n, \mu}(1)=\sum_{k=0}^n\underset{t=k}{\rm res}\,\left(
F(n,t)\, e^{i\pi t}\cdot\Bigl(\frac{\pi}{\sin\pi t}\Bigr)^{\mu+1}\right)
$$
and using the expansions (\ref{eq19})--(\ref{eq21}) we obtain
\begin{equation}
\begin{split}
{\rm Re}\, \widetilde{I}_{n,\mu}(1)&=\sum_{j=0}^{\mu/2}\frac{(-1)^j(2\pi)^{2j}
B_{2j}^{(\mu+1)}\bigl(\frac{\mu+1}{2}\bigr)}{(2j)!}\sum_{l=0}^{\mu/2-j}
\frac{(-1)^l \pi^{2l}}{(2l)!(\mu-2j-2l)!}
F_{n,\mu-2j-2l} \\
&=\frac{1}{\mu!}F_{n,\mu}+\sum_{\nu=1}^{\mu/2}\widetilde{d}_{2\nu,\mu} F_{n, \mu-2\nu},
\label{eq22}
\end{split}
\end{equation}
where $\widetilde{d}_{2\nu,\mu}$ are some real constants independent of $n$ and
$$
\widetilde{d}_{\mu,\mu}=(i\pi)^{\mu}\sum_{j=0}^{\mu/2}
\frac{4^j B_{2j}^{(\mu+1)}\bigl(\frac{\mu+1}{2}\bigr)}{(2j)!(\mu-2j)!}.
$$
Taking into account (\ref{eq17}) we can rewrite the last equality as
$$
\widetilde{d}_{\mu,\mu}=(i\pi)^{\mu}\sum_{j=0}^{\mu}\frac{2^j B_j^{(\mu+1)}\bigl(
\frac{\mu+1}{2}\bigr)}{j! (\mu-j)!}=\frac{(i\pi)^{\mu}}{\mu!}\sum_{j=0}^{\mu}\binom{\mu}{j}
B_j^{(\mu+1)}\Bigl(\frac{\mu+1}{2}\Bigr)2^j,
$$
which implies $\widetilde{d}_{\mu,\mu}=0,$ by Lemma \ref{l3}.
Now according to (\ref{eq22})   we have 
$$
{\rm Re}\, \widetilde{I}_{n,\mu}(1)=\frac{1}{\mu!}F_{n,\mu}+\sum_{\nu=1}^{\mu/2-1}\widetilde{d}_{2\nu,\mu} F_{n, \mu-2\nu},
$$
from which, by the induction hypothesis, the lemma follows.
\end{proof}
\begin{lemma} \label{l5}
Let $a,$ $\mu$ be integers satisfying $a\ge 2,$ $0\le\mu\le a-1$ and $u\in {\mathbb R}.$
Then for each $n=0,1,2,\ldots$ we have
$$
I_{n,\mu}(u)=\frac{1}{(2\pi i)^{a-\mu-1}}\sum_{k=0}^{a-\mu-1}
(-1)^k\binom{a-\mu-1}{k}I_{n,a-1}(a+u-\mu-1-2k).
$$
\end{lemma}
\begin{proof}
From (\ref{eq11}) and the reflection formula for the gamma function
$$
\Gamma(t-n) \Gamma(n-t+1)=\frac{(-1)^n \pi}{\sin\pi t}
$$
we have
\begin{equation*}
\begin{split}
I_{n,\mu}(u)&=\frac{n!^a}{2\pi i}\int_L\frac{e^{i\pi tu}}{\Gamma^{a-1}(t+1)\Gamma^a(n-t+1)}
\left(\frac{\pi}{\sin\pi t}\right)^{\mu+1}\,dt \\[3pt]
&=\frac{(-1)^{an} n!^a}{2\pi i}\int_L\frac{\Gamma^a(t-n)}{\Gamma^{a-1}(t+1)}\left(\frac{\sin\pi t}{\pi}
\right)^{a-\mu-1}e^{i\pi tu}\,dt.
\end{split}
\end{equation*}
Now replacing the function $\sin(\pi t)$ by $(e^{i\pi t}-e^{-i\pi t})/(2i)$
and opening parenthesis, by the binomial theorem, we get
\begin{equation*}
I_{n,\mu}(u)=\frac{(-1)^{an} n!^a}{(2\pi i)^{a-\mu}}\sum_{k=0}^{a-\mu-1}
(-1)^k\binom{a-\mu-1}{k}\int_L\frac{\Gamma^a(t-n)}{\Gamma^{a-1}(t+1)}\,
e^{i\pi t(u+a-\mu-1-2k)}\,dt,
\end{equation*}
which completes the proof of the lemma.
\end{proof}
\begin{lemma} \label{l6}
Let $a, \mu$ be integers satisfying  $a\ge 2,$ $0\le\mu\le a-1,$ and $u\in{\mathbb R}.$ Then
for each $n=0,1,2,\ldots$ we have
$$
I_{n, \mu}(-u)=\overline{I_{n, \mu}(u)},
$$
where the bar stands for complex conjugation.
\end{lemma}
\begin{proof}
Making the change of variable $t\mapsto\overline{t}$ in the integral $I_{n, \mu}(u)$
and using the equalities  $\Gamma(\bar{z})=\overline{\Gamma(z)}$ and
$\sin(\bar{z})=\overline{\sin(z)}$ we get
\begin{equation*}
\begin{split}
I_{n, \mu}(u)&=\frac{1}{2\pi i}\int_L F(n,\overline{t})
\left(\frac{\pi}{\sin(\pi \overline{t})}\right)^{\mu+1} e^{i\pi u\overline{t}}\,
d\overline{t}
=-\frac{1}{2\pi i}\int_L
F(n,\overline{t})
\left(\frac{\pi}{\sin(\pi \overline{t})}\right)^{\mu+1} e^{i\pi u\overline{t}}\,
dt \\[5pt]
&=\overline{\frac{1}{2\pi i}}\int_L
\overline{\left(\frac{\pi}{\sin(\pi t)}\right)^{\mu+1}}\overline{F(n,t)}\,
 e^{\overline{-i\pi ut}}\,dt=
\overline{\frac{1}{2\pi i}}\int_L\overline{F(n,t)\left(\frac{\pi}{\sin\pi t}\right)^{\mu+1}
e^{-i\pi ut}}\,dt \\[5pt]
&=\overline{I_{n, \mu}(-u)},
\end{split}
\end{equation*}
as required.
\end{proof}

\section{Asymptotics of the integral $I_{n,a-1}(u)$}

\vspace{0.3cm}

\begin{lemma} \label{l7}
Let $a, u\in{\mathbb Z},$ $a\ge 2,$ $|u|\le a,$ and $n$ be a sufficiently large
positive integer. Then all $a$ roots of the polynomial
$$
p_u(\tau)=e^{i\pi u}n(\tau-1)^a-\tau^{a-1}
$$
are given by the following asymptotic expansions:
$$
\tau_k(u)=1+\sum_{m=1}^{\infty}\frac{(2-m/a)_{m-1}}{m!}\cdot\frac{e^{\frac{m(2\pi k-\pi u)i}{a}}}{n^{m/a}},
\qquad k=0,1,\ldots, a-1.
$$
\end{lemma}
\begin{proof}
First note that the polynomial $p_u(\tau)$ has not  real roots on $(-\infty,0].$ Indeed, if we suppose
that $\tau=-x,$ $x\ge 0,$ is such a root, then $p_u(-x)=0$ and
\begin{equation}
(-1)^u n(x+1)^a+x^{a-1}=0.
\label{eq26}
\end{equation}
On the other hand, we have that the left-hand side of (\ref{eq26}) is positive (negative)
if $u$ is even (odd), which is a contradiction.  Therefore, we can consider the equation
\begin{equation}
n(\tau-1)^ae^{i\pi u}=\tau^{a-1}
\label{eq27}
\end{equation}
in the complex $\tau$-plane with cut along the ray $(-\infty; 0].$ Then it is easily seen
that the equation (\ref{eq27})  is equivalent to $a$ relations of the form
\begin{equation}
n^{\frac{1}{a}} e^{\frac{\pi u-2\pi k}{a}i} (\tau-1)=\tau^{1-\frac{1}{a}},
\qquad k=0,1,\ldots,a-1,
\label{eq23}
\end{equation}
where $\tau^{1-\frac{1}{a}}=e^{(1-\frac{1}{a})\log\tau}$ and  logarithm
is defined by its principal branch.  From (\ref{eq23})
we get immediately 
all $a$ roots of the polynomial $p_u(\tau):$
\begin{equation}
\tau_k(u)=1+\frac{\exp(\frac{2\pi k-\pi u}{a}i)}{n^{1/a}}
+
\frac{a-1}{a}\cdot\frac{\exp(\frac{2(2\pi k-\pi u)}{a}i)}{n^{2/a}}+O(n^{-3/a}), \quad
k=0,1,\ldots,a-1.
\label{eq24}
\end{equation}
It is possible to find a complete asymptotic expansion in (\ref{eq24}) if we apply
the Lagrange inversion formula (see \cite[\S 2.2]{debru}) to the equation (\ref{eq23}).
Indeed, substituting $\tau-1=z$ and rewriting (\ref{eq23}) as
\begin{equation}
\frac{z}{(z+1)^{1-\frac{1}{a}}}=n^{-\frac{1}{a}}\cdot e^{\frac{2\pi k-\pi u}{a}i},
\label{eq25}
\end{equation}
we get that there exist positive numbers $\rho_1$ and $\rho_2$ such that for
$n^{-1/a}<\rho_1$ the equation (\ref{eq25}) has just one solution $z$ in the domain
$|z|<\rho_2,$ that is
\begin{equation}
z=\sum_{m=1}^{\infty}\frac{c_m}{n^{m/a}}\, e^{\frac{m(2\pi k-\pi u)}{a}i},
\label{eq28}
\end{equation}
where the coefficients $c_m$ are given by the formula
\begin{equation}
c_m=\frac{1}{m!}\left.\left(\Bigl(\frac{d}{dz}\Bigr)^{m-1}(z+1)^{m(1-1/a)}\right)\right|_{z=0}.
\label{eq29}
\end{equation}
If we suppose that there is another solution of (\ref{eq25}) with $|z|\ge \rho_2,$
then since the function $\frac{|z|}{(|z|+1)^{1-1/a}}$
increases from $0$ to $+\infty$ as $|z|$ increases from $0$ to $+\infty,$ we get
$$
\frac{\rho_2}{(\rho_2+1)^{1-1/a}}\le
\frac{|z|}{(|z|+1)^{1-1/a}}\le
\frac{|z|}{|z+1|^{1-1/a}}=n^{-1/a},
$$
which is impossible if $n$ is sufficiently large.
Hence formulas (\ref{eq28}), (\ref{eq29}) give the complete asymptotic expansion for $\tau_k(u),$
$$
\tau_k(u)=1+\sum_{m=1}^{\infty}\frac{(2-m/a)_{m-1}}{m!\, n^{m/a}}\,
e^{\frac{m(2\pi k-\pi u)}{a}i}, \qquad k=0,1,\ldots, a-1,
$$
and the lemma is proved.
\end{proof}
\begin{lemma} \label{l8}
Let $a, u\in {\mathbb Z,}$ $a\ge 2,$ $|u|\le a.$ Then the following formula
holds  for the integral $I_{n,a-1}(u)$ as $n\to\infty:$
\begin{equation}
I_{n,a-1}(u)=\frac{(-1)^u(2\pi)^{\frac{a-1}{2}}}{in^{\frac{a-3}{2}}}\int_{L_1}
g(\tau)e^{nf(\tau)}\,d\tau\cdot (1+O(n^{1/a-1})),
\label{eq40}
\end{equation}
where
\begin{equation}
f(\tau)=a(\tau-1)\log(\tau-1)-(a-1)\tau\log\tau-\tau+\tau\log n+i\pi\tau u,
\label{eq41}
\end{equation}
\begin{equation}
\quad g(\tau)=\frac{(\tau-1)^{\frac{a}{2}}}{\tau^{\frac{3(a-1)}{2}}},
\label{eq42}
\end{equation}
and $L_1$ is a loop beginning and ending at $-\infty,$ encircling the points $1,0,-1,\ldots$
once in the positive direction and intersecting the real axis at the point $\tau_0:=\tau_0(0)$
if $u\ne 0,$ and $L_1$ is a vertical line going from $\tau_0-i\infty$ to $\tau_0+i\infty$ if $u=0.$
\end{lemma}
\begin{proof}
To make the many-valued  functions $f(\tau)$ and $g(\tau)$ definite, we consider them in the $\tau$-plane
with cut along the ray $(-\infty,1]$ fixing the branches of logarithms that take real values
on the interval $(1,+\infty)$ of the real axis.

Note that the contour of integration $L$
in the integral $I_{n, a-1}(u)$  defined in (\ref{eq11}) is an arbitrary loop
beginning and ending at $-\infty$ and
encircling the points $n, n-1, n-2, \ldots$ once in the positive direction
or an arbitrary vertical line ${\rm Re}\,\tau=c$ with $c>n.$
Now let us denote by $\tau_0:=\tau_0(0)=1+\frac{1}{n^{1/a}}+O(n^{-2/a})$ one of the real
roots of the polynomial $p_0(\tau)$ (see Lemma \ref{l7}) and suppose that the contour $L$
(loop or vertical line) intersects the real axis at the point $n\tau_0+1$ and for any point $t$ of $L,$
$|t|\ge n\tau_0+1.$

Then  the asymptotic expansion of the Gamma function for large $|z|$ (see \cite[\S 2.11]{lu}),
\begin{equation}
\log\Gamma(z)=\left(z-\frac{1}{2}\right)\log z-z+\frac{1}{2}\log(2\pi)+O(|z|^{-1}),
\label{gamma}
\end{equation}
where $|\arg z|\le \pi-\varepsilon,$
$\varepsilon>0$ and the constant in $O$ is independent of $z,$
implies
that the following
formula holds for the integrand of $I_{n, a-1}(u)$ on the contour $L:$
\begin{equation*}
\begin{split}
&\frac{\Gamma^a(t-n)}{\Gamma^{a-1}(t+1)}e^{i\pi tu}=\exp\{
a\log\Gamma(t-n)-(a-1)\log\Gamma(t+1)+i\pi tu\}\\
&=\exp\left\{a\Bigl(t-n-\frac{1}{2}\Bigr)\log(t-n)-(a-1)\Bigl(t+\frac{1}{2}\Bigr)\log(t+1)
+i\pi tu \right.\\
&+an-t+a-1\left.+\log\sqrt{2\pi}+O(n^{1/a-1})\right\}.
\end{split}
\end{equation*}
The change of variable $t=n\tau+1$ yields
\begin{equation*}
\begin{split}
\frac{\Gamma^a(t-n)}{\Gamma^{a-1}(t+1)}e^{i\pi tu}&=(-1)^u\sqrt{2\pi}\,
e^{an-an\log n-(a-3/2)\log n}\frac{(\tau-1)^{\frac{a}{2}}}{\tau^{\frac{3(a-1)}{2}}}\,
e^{nf(\tau)} (1+O(n^{1/a-1}))\\
&=\frac{(-1)^u(2\pi)^{\frac{a+1}{2}}}{n^{\frac{a-3}{2}}n!^a}\,g(\tau)\,e^{nf(\tau)} (1+O(n^{1/a-1})),
\end{split}
\end{equation*}
where the functions $f(\tau), g(\tau)$ are defined in (\ref{eq41}), (\ref{eq42}).
This completes the proof of the lemma.
\end{proof}
The next lemma is devoted to  the calculation of the asymptotics of $I_{n, a-1}(u)$ using the saddle-point
method.
\begin{lemma} \label{l9}
Let $a, u\in{\mathbb Z}, a\ge 2, |u|\le a.$ Then
the asymptotic behavior of the integral $I_{n, a-1}(u)$ as $n\to\infty$ is given by the formula:
$$
I_{n,a-1}(u)=\frac{(-1)^{nu}(2\pi)^{\frac{a-1}{2}}\,e^{\frac{\pi u(a-1)i}{2a}}}{\sqrt{a}}
  \frac{n!}{n^{\frac{a}{2}+\frac{1}{2a}}}\,
\exp\!\left(\sum_{m=1}^ab_m(a)\,e^{\frac{-\pi mui}{a}}n^{1-\frac{m}{a}}\right)
(1+O(n^{-1/a})),
$$
where  the coefficients $b_m(a)$ are rational numbers that can be found explicitly from the expansion
{\rm (\ref{bma})}.
\end{lemma}
\begin{proof}
It is sufficient to prove the lemma only for $u\le 0.$ The case when $u>0$ can be reduced to this case
by Lemma \ref{l6}.
First we determine the saddle points of the integrand (\ref{eq40}), which are the zeros of the derivative
of the function $f(\tau).$ It is easy to see that the zeros of the derivative
$$
f'(\tau)=a\log(\tau-1)-(a-1)\log \tau+\log n+i\pi u
$$
are simultaneously roots of the polynomial $p_u(\tau)$ defined in Lemma \ref{l6}.
For the root $\tau_k(u),$ $k=0,1,\ldots, a-1,$ by Lemma \ref{l6}, we have
$$
\tau_k(u)-1=\frac{e^{\frac{2\pi k-\pi u}{a}i}}{n^{1/a}}\left(1+
\frac{a-1}{a}\cdot\frac{e^{\frac{2\pi k-\pi u}{a}i}}{n^{1/a}}+O(n^{-2/a})\right).
$$
Since
$$
-\frac{\pi u}{a}\le\frac{2\pi k-\pi u}{a}\le -\frac{\pi u}{a}+2\pi -\frac{2\pi}{a},
$$
it is easily seen that $\tau_0(u)$ is exactly one zero of the derivative $f'(\tau)$ among
all the roots (\ref{eq24}). Therefore, $\tau_0(u)$ is the only saddle point of the function
$e^{nf(\tau)}.$

By a similar argument as in the proof of Lemma \ref{l6}, we can find complete asymptotic expansions
for the points $\tau$ satisfying the equation ${\rm Re} f'(\tau)=0,$ which is equivalent
to
\begin{equation}
\frac{|z|}{|z+1|^{1-1/a}}=n^{-1/a}, \qquad\mbox{where}\quad z=\tau-1.
\label{eq46}
\end{equation}
Application of the Lagrange inversion formula to the equation
$$
\frac{z}{(z+1)^{1-1/a}}=w
$$
shows that there are positive constants $\delta_1$ and $\delta_2$ such that for
$|w|<\delta_2$ there is only one solution $z$ satisfying $|z|<\delta_1$ and this solution
is an analytic function of $w:$
$$
z=\sum_{m=1}^{\infty}\frac{(2-m/a)_{m-1}}{m!} w^m
$$
(actually, by the asymptotic formula (\ref{gamma}) for the Gamma function, the series converges
if $|w|<a/(a-1)^{1-1/a}.$)
So it is clear that if $n>\delta_2^a,$ there is one and only one solution of (\ref{eq46})
in the circle $|\tau-1|<\delta_1$ given by
\begin{equation}
\tau=1+\sum_{m=1}^{\infty}\frac{(2-m/a)_{m-1}}{m!} \frac{e^{im\varphi}}{n^{m/a}}, \qquad
-\pi<\varphi\le\pi.
\label{eq47}
\end{equation}
On the other hand, if we suppose that there is another solution with $|z|\ge \delta_1,$
then by the same argument as in the proof of Lemma \ref{l6}, we get a contradiction
if $n$ is sufficiently large. So if $n$ is large enough, the continuously
differentiable curve
$$
{\rm Re}f'(\tau)=\log \frac{|\tau-1|^an}{|\tau|^{a-1}}=0
$$
represents approximately a small circle with center at $\tau=1$ and radius $n^{-1/a}.$ The points of this curve
are given by their complete asymptotic expansion (\ref{eq47}).
Moreover, we have ${\rm Re} f'(\tau)<0$ inside of this curve and ${\rm Re} f'(\tau)>0$ outside of it.
It is easily seen that ${\rm Re} f'(\tau_k(u))=0$ for any integers $u$ and $k.$

To apply the saddle-point method  to evaluation of the integral $I_{n, a-1}(u),$
we need to choose a new contour of integration $L_2$ passing through the saddle point $\tau_0(u).$
If $u=0$  we use the original contour ${\rm Re}\,\tau=\tau_0.$
 Then for $\tau=\tau_0+iy,$ $-\infty<y<+\infty$ by the
Cauchy-Riemann conditions, we have
\begin{equation*}
\frac{d}{dy}\,{\rm Re} f(\tau_0+iy)=-{\rm Im}\frac{d}{d\tau} f(\tau_0+iy)=
-a\arg(\tau-1)+(a-1)\arg\tau.
\end{equation*}
Since ${\rm Re}\, \tau_0>1,$ for $y<0$ we get
$$
-\frac{\pi}{2}<\arg(\tau-1)<\arg\tau<0
$$ and therefore,
\begin{equation}
\frac{d}{dy}\, {\rm Re} f(\tau_0+iy)=(a-1)(\arg \tau-\arg(\tau-1))-\arg(\tau-1)>0.
\label{eq48}
\end{equation}
This implies that ${\rm Re} f(\tau_0+iy)$ strictly increases as $y$ increases from $-\infty$ to $0.$
If $y>0$ then we have
$$
0<\arg\tau<\arg(\tau-1)<\frac{\pi}{2}
$$
and hence
$$
\frac{d}{dy}\, {\rm Re} f(\tau_0+iy)=(a-1)(\arg\tau-\arg(\tau-1))-\arg(\tau-1)<0.
$$
Therefore, the function ${\rm Re} f(\tau_0+iy)$ strictly decreases as $y$ increases from $0$ to
$+\infty.$ This proves that ${\rm Re} f(\tau_0+iy)$ attains its maximum on $L_1$
at the unique point $\tau_0$ and we can apply the saddle-point method to calculate
the asymptotics of $I_{n,a-1}(0)$ (see (\ref{eq60}) below).

If $u\ne 0$ we define $L_2$ as a contour consisting of three parts:

$(i)\,$ the half-line $\tau=\tau_0+iy,$ $-\infty<y\le 0;$

$(ii)\,$ the segment $\tau=\tau_0+re^{i\varphi},$ $\varphi=\arg(\tau_0(u)-\tau_0),$
$0\le r\le |\tau_0(u)-\tau_0|,$  connecting the points $\tau_0$ and $\tau_0(u);$

$(iii)\,$ the half-line $\tau=-x+i\,{\rm Im}\, \tau_0(u),$ $-{\rm Re}\,\tau_0(u)\le x<+\infty.$

If $u=-a$ both parts $(ii)$ and $(iii)$ give a ray going from $\tau_0$ to $-\infty$ along the
upper bank of the cut $[-\infty,1].$

Now  show that we can replace the contour $L_1$ in the integral $I_{n,a-1}(u)$ by the contour $L_2,$
i.e., show that
\begin{equation}
\int_{L_1}e^{nf(\tau)} g(\tau)\,dt=\int_{L_2} e^{nf(\tau)} g(\tau)\, dt.
\label{eq52}
\end{equation}
For this purpose,  consider the circle $|\tau|=N,$ where $N$ is a sufficiently large integer.
Suppose that $L_1^{-},$ $L_2^{-}$ ($L_1^{+},$ $L_2^{+}$) are points of intersections of this circle
with the contours $L_1,$ $L_2$ in the lower half-plane (upper half-plane), respectively. Then
to prove (\ref{eq52}), it is sufficient to show that
\begin{equation}
\int_{L_1^{-}L_2^{-}} e^{nf(\tau)} g(\tau)\, dt\to 0, \qquad
 \int_{L_2^{+}L_1^{+}} e^{nf(\tau)} g(\tau)\, dt\to 0
\label{eq53}
\end{equation}
as $N\to\infty,$ where $L_1^{-}L_2^{-},$ $L_2^{+}L_1^{+}$ are arcs of the circle of radius $N$
with center at the origin. On the arcs  $L_1^{-}L_2^{-}$ and $L_2^{+}L_1^{+}$ of the circle $\tau=Ne^{i\varphi}$
for $N$ sufficiently large we have the inequalities $-\pi<\varphi<-\pi/4$ and $3\pi/4<\varphi\le\pi,$
respectively (the value $\varphi=\pi$ corresponds to the upper bank of the cut $(-\infty,1]$). By Taylor's
formula, we obtain
$$
\log(\tau-1)=\log(Ne^{i\varphi}-1)=\log N+i\varphi-\frac{e^{-i\varphi}}{N}+O(N^{-2}),
$$
where the constant in $O(N^{-2})$ is absolute.  Substituting this expansion in
(\ref{eq41}), we get
$$
{\rm Re}\,f(Ne^{i\varphi})=
N\log N\cos(\varphi)+N\log\frac{n}{e}\cos(\varphi)-N(\varphi+\pi u)\sin(\varphi)-a\log N
-a+O(N^{-1}).
$$
Note that on the arc $L_2^{-}L_1^{-}$ we have $-1\le\cos(\varphi)\le\tau_0/\sqrt{N^2+\tau_0^2},$
$\sin(\varphi)\le 0$  and $\varphi+\pi u<0.$ This yields
\begin{equation}
{\rm Re}\,f(\tau)<\tau_0\log N+\tau_0\log\frac{n}{e}-a\log N-a+O(N^{-1})=
(\tau_0-a)\log N+O(1).
\label{eq54}
\end{equation}
Similarly, on the arc $L_1^{+}L_2^{+}$ we have $-1\le\cos(\varphi)\le -\sqrt{2}/2,$ $\sin(\varphi)\ge 0,$
and therefore,
\begin{equation}
{\rm Re}\, f(\tau)<-\frac{\sqrt{2}}{2}N\log N+O(N).
\label{eq55}
\end{equation}
For the function $g(\tau)$ on the arcs $L_1^{-}L_2^{-}$ and
$L_1^{+}L_2^{+},$ we have the following trivial estimate:
\begin{equation}
|g(\tau)|=O(N^{3/2-a}) \qquad\mbox{as}\quad N\to\infty.
\label{eq56}
\end{equation}
Since the length of each of the arcs $L_1^{-}L_2^{-}$ and $L_1^{+}L_2^{+}$ does not exceed $\pi N,$
estimates (\ref{eq54})--(\ref{eq56}) imply that the integrals in (\ref{eq53}) are of orders
$O(N^{\tau_0+1/2-2a})$ and $O(N^{1/2-a-\sqrt{2}/2N}),$ respectively. Hence, the limiting relations in
(\ref{eq53}) hold, and we can replace the contour of integration $L_1$ in the integral $I_{n,a-1}(u)$
by the new contour $L_2.$

Now we show that $\tau_0(u)$ is the unique maximum point of ${\rm Re}\,f(\tau)$ on the contour $L_2.$
Since $u<0,$ from (\ref{eq48}) it follows that $\frac{d}{dy}\,{\rm Re} f(\tau_0+iy)$  is positive when
$y<0$ and therefore, ${\rm Re} f(\tau)$ monotonically increases on the half-line $\tau=\tau_0+iy,$
$-\infty<y\le 0.$ Similarly, for the half-line $\tau=-x+i\,{\rm Im}\, \tau_0(u)$ we have
$$
\frac{d}{dx}\, {\rm Re} f(-x+i\,{\rm Im}\,\tau_0(u))=-{\rm Re}\, \frac{d}{d\tau}
f(-x+i\,{\rm Im}\,\tau_0(u))<0.
$$
This shows that the function ${\rm Re}\,f(\tau)$ monotonically decreases on the half-line
going from $\tau_0(u)$ to $-\infty+i\,{\rm Im}\, \tau_0(u).$

Now consider the segment $\tau=\tau_0+r e^{i\varphi}$ defined in $(ii).$ The derivative of
${\rm Re}\, f(\tau)$ on this part of $L_2$ is given by the formula:
\begin{equation}
\frac{d}{dr}\,{\rm Re}\,f(\tau)={\rm Re}\,\frac{df(\tau)}{dr}={\rm Re}\,
\left(\frac{df(\tau)}{d\tau}\cdot\frac{d\tau}{dr}\right)=
{\rm Re}\,f'(\tau) \cos\varphi-{\rm Im}\,f'(\tau) \sin\varphi.
\label{eq49}
\end{equation}
Note that since ${\rm Re}\,\tau_0(u)<{\rm Re}\,\tau_0,$ it follows that $\pi/2<\varphi<\pi,$
and therefore the product ${\rm Re}\,f'(\tau)\cdot \cos\varphi$ is positive on the segment.
Let us investigate the behavior of ${\rm Im}\,f'(\tau)$ on this part of the contour. Note that
${\rm Im}\,f'(\tau_0)=\pi u<0$ and ${\rm Im}\,f'(\tau_0(u))=0.$ We show  that ${\rm Im}\,f'(\tau)$
monotonically increases on our segment from $\pi u$ to $0.$ To see this, we consider the derivative
of ${\rm Im}\,f'(\tau).$ By the Cauchy-Riemann conditions, for $r>0,$ we have
$$
\frac{d}{dr}\,{\rm Im}\,f'(\tau)={\rm Im}\,(f''(\tau)\cdot e^{i\varphi})=\frac{1}{r}\,
{\rm Im}\,(f''(\tau)\cdot(\tau-\tau_0)).
$$
Taking into account that
$$
f''(\tau)=\frac{a}{\tau-1}-\frac{a-1}{\tau}
$$
we obtain
$$
\frac{d}{dr}\,{\rm Im}\,f'(\tau)=\frac{1}{r}\, {\rm Im}\,
\left(\frac{a(\tau-\tau_0)}{\tau-1}-\frac{(a-1)(\tau-\tau_0)}{\tau}\right)=
\frac{{\rm Im}\,\tau}{r}\left(\frac{a(\tau_0-1)}{|\tau-1|^2}-
\frac{(a-1)\tau_0}{|\tau|^2}\right).
$$
Since $\tau$ is in the upper half-plane, to show that $\frac{d}{dr}\,{\rm Im}\,f'(\tau)$ is positive
on the segment, it is sufficient  to show that the quantity in parenthesis, which we denote by $B,$
is positive. For $\tau=\tau_0+re^{i\varphi},$ we have
$$
B:=\frac{a(\tau_0-1)}{|\tau-1|^2}-
\frac{(a-1)\tau_0}{|\tau|^2}=\frac{h(r)}{|\tau-1|^2 |\tau|^2},
$$
where
$$
h(r)=(\tau_0-a)r^2+2r\tau_0(\tau_0-1)\cos\varphi+\tau_0(\tau_0-1)(a+\tau_0-1).
$$
Since $\tau_0<a,$ it follows easily that the quadratic polynomial $h(r)$ has two real roots
$r_1,$  $r_2$ such that $r_1<0<r_2$ and $h(r)$ is positive on $(r_1, r_2)$ and negative
on $(-\infty, r_1),$ $(r_2, +\infty).$
Now we show that the point $r(u):=|\tau_0(u)-\tau_0|$ belongs to the interval $(0, r_2).$
Indeed, taking into account that
$$
|\tau_0-1|=\frac{1}{n^{1/a}}+O(n^{-2/a})\qquad\mbox{and}\qquad
|\tau_0(u)-1|=\frac{1}{n^{1/a}}+O(n^{-2/a})
$$
for $n$ sufficiently large, we have $|\tau_0(u)-1|\le 2|\tau_0-1|$ and $\tau_0<9/8,$
and therefore,
$$
r(u)=|\tau_0(u)-\tau_0|\le |\tau_0(u)-1|+|\tau_0-1|\le 3(\tau_0-1).
$$
This implies that
\begin{equation*}
\begin{split}
h(r(u))&\ge 9(\tau_0-a) (\tau_0-1)^2-6\tau_0(\tau_0-1)^2+\tau_0(\tau_0-1)(a+\tau_0-1) \\
&=
(\tau_0-1)(a(9-8\tau_0)+4\tau_0(\tau_0-1))>0
\end{split}
\end{equation*}
and hence  $r(u)\in (0, r_2).$ This proves that $\frac{d}{dr}\,{\rm Im}\,f'(\tau)$ is positive
on the segment $[0, r(u)]$ and therefore, ${\rm Im}\,f'(\tau)$ increases on it from $\pi u$ to $0,$
i.~e., ${\rm Im}\,f'(\tau)<0$ on $[0, r(u)).$ Now by (\ref{eq49}), we get $\frac{d}{dr}\,{\rm Re}\,f(\tau)>0$
on the segment connecting $\tau_0$ and $\tau_0(u)$ and therefore, ${\rm Re}\,f(\tau)$
is monotonically increasing on it.
Hence we showed that $\tau_0(u)$ is the unique maximum point on the whole contour of integration $L_2$
and we can apply the saddle-point method (see \cite[Ch.6, Th.3.1]{fe}) to estimate the integral $I_{n, a-1}(u).$

Let us recall that when $u=-a$ the contour $L_2$ consists of the two parts: vertical half-line
$\tau=\tau_0+iy,$ $-\infty<y\le 0,$ and horizontal half-line $\tau=-x+i 0,$ $x\ge {\rm Re}\,\tau_0(a),$
going along the upper bank of the cut $(-\infty, 1].$
Let us notice that in this case  the function $g(\tau)$
has singularity at the point $\tau=0$ belonging to the contour $L_2.$ So in order
to apply the saddle-point method in this case,  we should change slightly
the contour $L_2$  to avoid this singularity. It is easy to see that instead of the horizontal line
$\tau=-x+i0,$ $x\ge -1/2,$ we can take a contour consisting of
 the semi-circle $\tau=1/2e^{i\varphi},$ $0\le\varphi\le\pi,$
and the ray $\tau=-x+i0,$ $x\ge 1/2.$ Then we get
$$
\int_{L_2}e^{nf(\tau)} g(\tau)\,dt=\int_{\tau_0-i\infty}^{\tau_0+i0}e^{nf(\tau)} g(\tau)\,dt
+\int_{\tau_0+i0}^{1/2+i0}e^{nf(\tau)} g(\tau)\,dt+O(C^nn^{1/2}),
$$
where $C=C(a)$ is some positive constant independent of $n.$ So we can apply the saddle-point
method in this case to the contour $L_{2,a}:=\{\tau_0+iy, -\infty<y\le 0\}\cup\{x+i0,
1/2\le x\le \tau_0\}$ and as we can see from (\ref{eq57}) below the quantity $O(C^nn^{1/2})$
does not influence on the contribution of the saddle point.

Finally, applying the saddle-point method  we have
\begin{equation}
I_{n,a-1}(u)
=\frac{(-1)^u(2\pi)^{\frac{a}{2}}}{i n^{\frac{a-2}{2}}} \,
e^{\frac{\pi i}{2}-\frac{i}{2}\arg f''(\tau_0(u))}|f''(\tau_0(u))|^{-\frac{1}{2}}
g(\tau_0(u))e^{nf(\tau_0(u))}(1+O(n^{-\frac{1}{a}})).
\label{eq60}
\end{equation}
In order to find the contribution of the saddle point, we evaluate
\begin{equation*}
\begin{split}
f(\tau_0(u))&=-a\log(\tau_0(u)-1)-\tau_0(u)=\pi ui+\log n-1-
\sum_{m=1}^a\frac{(2-\frac{m}{a})_{m-1}}{m!} \, \frac{e^{\frac{-\pi mui}{a}}}{n^{m/a}} \\
&-a\log\Bigl(1+\sum_{m=1}^a\frac{(2-\frac{m+1}{a})_m}{(m+1)!}\,\frac{e^{\frac{-\pi mui}{a}}}{n^{m/a}}\Bigr)
+O(n^{-1-1/a}).
\end{split}
\end{equation*}
Now expanding the logarithm in powers of $e^{\frac{-\pi ui}{a}} n^{-1/a}$ we get
\begin{equation}
\begin{split}
-a\log\Bigl(1&+\sum_{m=1}^a\frac{(2-\frac{m+1}{a})_m}{(m+1)!}\,\frac{e^{\frac{-\pi mui}{a}}}{n^{m/a}}\Bigr)
-\sum_{m=1}^a\frac{(2-\frac{m}{a})_{m-1}}{m!} \, \frac{e^{\frac{-\pi mui}{a}}}{n^{m/a}}\\
&=\sum_{m=1}^a b_m(a) \frac{e^{\frac{-\pi mui}{a}}}{n^{m/a}}+ O(n^{-1-1/a}),
\label{eq61}
\end{split}
\end{equation}
where $b_m(a)$ are rational coefficients depending only on $a,$ in particular, $b_1(a)=-a,$
$b_2(a)=(1-a)/2,$ $b_3(a)=(1-a)(2a-3)/(6a).$  From (\ref{eq61}) it follows that
\begin{equation}
e^{nf(\tau_0(u))}=\frac{(-1)^{un} n!}{\sqrt{2\pi n}}\, \exp\left(\sum_{m=1}^ab_m(a)\,e^{\frac{-\pi mui}{a}}
n^{1-\frac{m}{a}}\right) (1+O(n^{-1/a})).
\label{eq57}
\end{equation}
Taking into account that
$$
g(\tau_0(u))=e^{\frac{-\pi ui}{2}}\, n^{-1/2} (1+O(n^{-1/a})),
\qquad
f''(\tau_0(u))=a e^{\frac{\pi ui}{a}}\,n^{1/a} (1+O(n^{-1/a})),
$$
we obtain the asymptotic behavior of $I_{n, a-1}(u):$
$$
I_{n, a-1}(u)=\frac{(-1)^{nu}(2\pi)^{\frac{a-1}{2}}\,e^{\frac{\pi u(a-1)i}{2a}}}{\sqrt{a}}
  \frac{n!}{n^{\frac{a}{2}+\frac{1}{2a}}}\,
\exp\!\left(\sum_{m=1}^ab_m(a)\,e^{\frac{-\pi mui}{a}}n^{1-\frac{m}{a}}\right)
(1+O(n^{-1/a})),
$$
and the lemma is proved.
\end{proof}

\section{Proof of Theorem 1}

\begin{lemma} \label{l10}
Let $a,$ $\mu$ be positive integers satisfying $a\ge 2,$ $0\le\mu\le a-1.$
Then there exist positive constants $\lambda_0=\lambda_0(a),$ $\lambda_1=\lambda_1(a)$ such that for every
positive integer $n,$
\begin{equation}
|I_{n,\mu}(0)|\le\frac{\lambda_0 n!}{n^{\frac{a}{2}+\frac{1}{2a}}} \exp\left(
\sum_{m=1}^{a-1}b_m(a)\cos\Bigl(\frac{\pi m(a-\mu-1)}{a}\Bigr) n^{1-\frac{m}{a}}\right)
\label{eq70}
\end{equation}
and for infinitely many positive integers $n$ the similar lower bound holds:
\begin{equation}
|I_{n,\mu}(0)|\ge\frac{\lambda_1 n!}{n^{\frac{a}{2}+\frac{1}{2a}}} \exp\left(
\sum_{m=1}^{a-1}b_m(a)\cos\Bigl(\frac{\pi m(a-\mu-1)}{a}\Bigr) n^{1-\frac{m}{a}}\right)
\label{eq71}
\end{equation}
Moreover, the following asymptotic formula:
\begin{equation*}
|I_{n,\mu}(1)|=\frac{n!}{\sqrt{a}\,(2\pi)^{\frac{a-1}{2}-\mu}\,
n^{\frac{a}{2}+\frac{1}{2a}}}\, \exp\left(
\sum_{m=1}^{a-1}b_m(a)\cos\Bigl(\frac{\pi m(a-\mu)}{a}\Bigr) n^{1-\frac{m}{a}}\right)
(1+O(n^{-1/a}))
\end{equation*}
as $n\to\infty$ takes place.
\end{lemma}
\begin{proof}
From Lemma \ref{l5} we have
\begin{equation}
I_{n,\mu}(1)=\frac{1}{(2\pi i)^{a-\mu-1}} \sum_{k=0}^{a-\mu-1}(-1)^k \binom{a-\mu-1}{k}
I_{n,a-1}(a-\mu-2k).
\label{eq72}
\end{equation}
By Lemmas \ref{l9}, \ref{l6}, it follows that the sum on the right of (\ref{eq72})
contains exactly one term with dominant asymptotics, that is $I_{n,a-1}(a-\mu).$
Therefore, we have
\begin{equation*}
\begin{split}
& \qquad\qquad\qquad\qquad  |I_{n,\mu}(1)|\sim \frac{1}{(2\pi)^{a-\mu-1}} |I_{n,a-1}(a-\mu)| \\
&=
\frac{n!}{\sqrt{a}\,(2\pi)^{\frac{a-1}{2}-\mu}\,
n^{\frac{a}{2}+\frac{1}{2a}}}\, \exp\left(
\sum_{m=1}^{a-1}b_m(a)\cos\Bigl(\frac{\pi m(a-\mu)}{a}\Bigr) n^{1-\frac{m}{a}}\right)
(1+O(n^{-1/a})),
\end{split}
\end{equation*}
as required. Similarly, for the integral $I_{n,\mu}(0)$ by Lemma \ref{l5}, we get
\begin{equation}
I_{n,\mu}(0)=\frac{1}{(2\pi i)^{a-\mu-1}} \sum_{k=0}^{a-\mu-1}(-1)^k \binom{a-\mu-1}{k}
I_{n,a-1}(a-\mu-1-2k).
\label{eq73}
\end{equation}
From Lemma \ref{l6} it follows that if $a-\mu$ is odd, then the quantity $I_{n,\mu}(0)$
represents a linear combination of $I_{n,a-1}(0)$ and sums of complex conjugates
$I_{n,a-1}(u)$ and $\overline{I}_{n,a-1}(u)$ for $u=2,4,\ldots, a-\mu-1.$
Similarly, if $a-\mu$ is even, then  $I_{n,\mu}(0)$ is equal to a linear combination
of differences of complex conjugates $I_{n,a-1}(u)$ and $\overline{I}_{n,a-1}(u)$
for $u=1,3,\ldots, a-\mu-1.$ Now by Lemma \ref{l9}, it is clear that the term
with dominant exponent on the right of (\ref{eq73}) is ${\rm Re}\, I_{n,a-1}(a-\mu-1)$
or ${\rm Im}\, I_{n,a-1}(a-\mu-1)$ depending on $a-\mu$ is odd or even.
From Lemma \ref{l6} we obtain
\begin{equation*}
\begin{split}
{\rm Re}\,I_{n,a-1}(a-\mu &-1)=\frac{(2\pi)^{\frac{a-1}{2}}}{\sqrt{a}}
\frac{n!}{n^{\frac{a}{2}+\frac{1}{2a}}}\,\cos(P(n)) \\
&\times\exp\left(\sum_{m=1}^ab_m(a)\cos\Bigl(
\frac{\pi m(a-\mu-1)}{a}\Bigr)n^{1-\frac{m}{a}}\right) (1+O(n^{-1/a}))
\end{split}
\end{equation*}
and
\begin{equation*}
\begin{split}
{\rm Im}\,I_{n,a-1}(a-\mu &-1)=\frac{(-1)^n(2\pi)^{\frac{a-1}{2}}}{\sqrt{a}}
\frac{n!}{n^{\frac{a}{2}+\frac{1}{2a}}}\, \sin(P(n)) \\
&\times\exp\left(\sum_{m=1}^ab_m(a)\cos\Bigl(
\frac{\pi m(a-\mu-1)}{a}\Bigr)n^{1-\frac{m}{a}}\right) (1+O(n^{-1/a})),
\end{split}
\end{equation*}
where
$$
P(n)=\frac{\pi(a-\mu-1)(a-1)}{2a}-\sum_{m=1}^ab_m(a)\sin\Bigl(
\frac{\pi m(a-\mu-1)}{a}\Bigr)n^{1-\frac{m}{a}}.
$$
Since the functions sine and cosine are bounded, we get immediately the required upper bound
(\ref{eq70}). On the other hand, by Weil's theorem, it is possible to show  (see Rivoal's
argument in the proof of Proposition 13 \cite{ri}) that each of the sequences
$\cos(P(n))$ and $\sin(P(n))$
is dense in the interval $[-1,1].$ This implies that  there are infinitely many $n$ such that the absolute
values of the cosine (sine) are not less than $1/2$ and hence the lower bound (\ref{eq71})
follows for infinitely many $n.$
\end{proof}
\begin{lemma} \label{l11}
Let $a\ge 2$ be an integer. Then there exists a positive constant $\lambda=\lambda(a)$
such that for every $\mu=0,1,\ldots, a-1,$ we have
$$
|F_{n,\mu}|\le\frac{\lambda\, n!}{n^{\frac{a}{2}+\frac{1}{2a}}}\,\exp\left(
\sum_{m=1}^{a-1}(-1)^m b_m(a) \cos\Bigl(\frac{2\pi m}{a}\Bigr) n^{1-\frac{m}{a}}\right).
$$
Moreover, the following asymptotic formula holds:
$$
q_n=\frac{n!}{\sqrt{a} (2\pi)^{\frac{a-1}{2}}\, n^{\frac{a}{2}+
\frac{1}{2a}}}\,\exp\left(
\sum_{m=1}^{a}(-1)^m b_m(a)  n^{1-\frac{m}{a}}\right)
(1+O(n^{-1/a}))
$$
as $n\to\infty.$
\end{lemma}
\begin{proof}
The proof follows easily from Lemmas \ref{l4}, \ref{l2}, \ref{l10} and the fact
that $q_n=F_{n,0}=I_{n,0}(1).$
\end{proof}
Now by Lemmas \ref{l1}, \ref{l11}, the theorem  follows.

\vspace{0.3cm}

{\bf\small Acknowledgements.} This research was in part supported by grants no.~89110024 (first author) and
no.~89110025 (second author) from School of Mathematics, Institute for Research in Fundamental Sciences (IPM).
This work was done during our summer visit in 2010 at the Abdus Salam International Centre for Theoretical Physics (ICTP),
Trieste, Italy. The authors thank the staff and personally the Head of the Mathematics Section of ICTP
Prof.~Ramadas Ramakrishnan for the hospitality and excellent working conditions. The first author is grateful to
the Commission on Development and Exchanges of the IMU for travel support.

\end{document}